\documentclass[a4paper,fleqn]{cas-dc}

\usepackage[numbers]{natbib}
\usepackage{amsthm}
\usepackage{url}
\usepackage{amssymb}
\usepackage{amsmath}
\usepackage{graphicx}

\newtheorem{theorem}{Theorem}
\newtheorem{problem}{Problem}

\usepackage{multirow}
\usepackage{algorithm}
\usepackage{algorithmic}
\usepackage{setspace}

\def\tsc#1{\csdef{#1}{\textsc{\lowercase{#1}}\xspace}}
\tsc{WGM}
\tsc{QE}
\tsc{EP}
\tsc{PMS}
\tsc{BEC}
\tsc{DE}

\begin{document}
\shortauthors{Junxiang Wang and Liang Zhao.}
\shorttitle{Nonconvex Generalization of Alternating Direction Method of Multipliers for Nonlinear Equality Constrained Problems}
\title [mode = title]{Nonconvex Generalization of Alternating Direction Method of Multipliers for Nonlinear Equality Constrained Problems}                      
\author[1]{Junxiang Wang}
\ead{jwan936@emory.edu}
\author[1]{Liang Zhao}
\ead{lzhao413@emory.edu}

\address[1]{Emory University, 201 Dowman Dr, Atlanta, GA, USA 30322}

\begin{abstract}
The classic Alternating Direction Method of Multipliers  (ADMM) is a popular framework to solve linear-equality constrained problems. In this paper, we extend the ADMM naturally to nonlinear equality-constrained problems, called neADMM. The difficulty of neADMM is to solve nonconvex subproblems. We provide globally optimal solutions to them in two important applications. Experiments on synthetic and real-world datasets demonstrate excellent performance and scalability of our proposed neADMM over existing state-of-the-start methods.
\end{abstract}

\maketitle
\section{Introduction}
\label{sec:introduction}
 \indent  There is a growing demand for efficient computational methods for analyzing high-dimensional large-scale data across
a wide variety of applications, including healthcare, finance, social media, astronomy, and e-commerce \cite{baniasadi2020transformation,wang2018multi,wang2018semi}. The classic Alternating Direction Method of Multipliers (ADMM) has received a significant amount of attention in the last few years. Its main advantage consists in the ability to split a complex problem into a series of simpler \textit{subproblems}, each of which is easy to solve \cite{boyd2011distributed}. While ADMM  focuses on optimization problems with linear equality constraints, many real-world problems require nonlinear constraints such as collaborative filtering \cite{candes2012exact}, 1-bit compressive sensing \cite{boufounos20081}, and mesh processing \cite{neumann2014compressed} and as yet, there lacks a discussion on how to apply ADMM to \emph{nonlinear equality-constrained} problems.\\
 \indent In this paper, we extend the ADMM  into the nonlinear equality-constrained problems in Section \ref{sec:neADMM}, called neADMM. Two important applications of neADMM are discussed in Section \ref{sec:application}, along with a consideration of ways to solve nonconvex subproblems.  Section \ref{sec:experiments} presents experiments conducted to show the effectiveness of the proposed neADMM on both synthetic and real-world datasets. We summarize this paper by Section \ref{sec:conclusions}.
\section{Nonlinear Equality-constrained ADMM}
\label{sec:neADMM}
We consider the following nonconvex problem with vector variables $x_1\in \mathbb{R}^{m_1}$ and $x_2\in \mathbb{R}^{m_2}$.
\begin{problem}
\label{prob:problem 2}
\begin{align*}
\nonumber &\min\nolimits_{x_1,x_2} \ F_1(x_1)+F_2(x_2) 
  \ s.t. \ f_1(x_1)+f_2(x_2)=0 
\end{align*}
\end{problem}
In Problem \ref{prob:problem 2}, $F_1(x_1)$ and $F_2(x_2)$  are proper continuous functions, and $f_1:\mathbb{R}^{m_1}\rightarrow \mathbb{R}^{d}$ and $f_2:\mathbb{R}^{m_2}\rightarrow \mathbb{R}^{d}$ can be nonlinear.
 We present the neADMM algorithm to solve Problem \ref{prob:problem 2}. According to the standard ADMM routine, we formulate the augmented Lagrangian as follows:
\begin{align}
 \nonumber &L_\rho(x_1,x_2,y)=F_1(x_1)+F_2(x_2)+y^T(f_1(x_1)+f_2(x_2))\\&+(\rho/2)\Vert f_1(x_1)+f_2(x_2)\Vert_2^2 \label{eq: augmented Lagrangian}
\end{align}
where $\rho>0$ is a penalty parameter and $y$ is a dual variable. neADMM aims to optimize the following two subproblems alternately:
\begin{align}
&x_1^{k+1} =\mathop {\arg \min }\nolimits_{{x_1}} L_\rho(x_1,x_2^k,y^k) \label{eq:subproblem 1}\\
&x_2^{k+1} =\mathop {\arg \min }\nolimits_{{x_2}} L_\rho(x^{k+1}_1,x_2,y^k) \label{eq:subproblem 2}    
\end{align}
Without loss of generality, we implicitly assume that there exist minima in Equations \eqref{eq:subproblem 1} and  \eqref{eq:subproblem 2}.\\
\indent The neADMM algorithm is presented in Algorithm \ref{algo:neADMM}. Specifically, Lines 3-4 update Equations \eqref{eq:subproblem 1} and \eqref{eq:subproblem 2}, Line 5 updates the dual variable $y$, and Lines 6 and 7 update the primal residual $r$ and the dual residual $s$, respectively.\\
\indent The main challenge of the neADMM framework is to solve nonconvex subproblems Equations \eqref{eq:subproblem 1} and \eqref{eq:subproblem 2}. While there is no general method to solve them exactly, for some specific forms, we have efficient solutions, which are discussed in the next section.
\begin{algorithm}  
\scriptsize
\caption{the neADMM Algorithm} 
\begin{algorithmic}[1] 
\STATE Initialize $x_1$ and $x_2$, $y$,  $\rho$, $k=0$.
\REPEAT
\STATE Update $x_1^{k+1}$ in Equation \eqref{eq:subproblem 1}.
\STATE Update $x_2^{k+1}$ in Equation \eqref{eq:subproblem 2}.
\STATE Update $y^{k+1}\leftarrow y^{k} + \rho(f_1(x_1^{k+1})+f_2(x_2^{k+1}))$
\STATE Update $r^{k+1}\leftarrow f_1(x^{k+1}_1)+f_2(x^{k+1}_2)$. $\#$ Calculate the primal residual.
\STATE Update $s^{k+1}\leftarrow \rho \partial f_1(x_1^{k+1})^T(f_2(x_2^{k+1})-f_2(x_2^{k}))$. $\#$ Calculate the dual residual.
\STATE $ k\leftarrow k+1$.
\UNTIL convergence.
\STATE Output $x_1$ and $x_2$.
\end{algorithmic}
\label{algo:neADMM}
\normalsize
\end{algorithm}

\section{Applications}
\label{sec:application}
\subsection{Optimization Problems with Spherical Constraints}
\indent The spherical constraint  is widely applied in 1-bit compressive sensing  \cite{boufounos20081} and mesh processing  \cite{neumann2014compressed}, which is formulated as follows:
\begin{problem}[Spherical Constrained Problem]
\label{prob:problem 4}
\begin{align*}
\min\nolimits_x \ \ell(x),\ \ \ \ s.t. \ \Vert x\Vert^2_2=1
\end{align*}
\end{problem}
\noindent where $\ell(\bullet)$ is a loss function.  We introduce an auxiliary variable $w$  and reformulate this problem as follows:
\begin{align*}
\min\nolimits_{x,w} \ell(x),\ \ \ \ s.t. \ \Vert w \Vert^2_2=1, w=x
\end{align*}
The augmented Lagrangian is formulated as follows according to Equation \eqref{eq: augmented Lagrangian}:
\begin{align*}
&L_\rho(x,w,y_1,y_2)\\&=\ell(x)+y^T_1(\Vert w\Vert^2_2-1)+(\rho/2)\Vert \Vert w\Vert^2_2-1\Vert^2_2\\&+y^T_2(w-x)+(\rho/2)\Vert w-x\Vert^2_2\\&\!=\!\ell(x)\!+\!(\rho/2)\Vert \Vert w\Vert^2_2-1+y_1/\rho\Vert^2_2\\&+(\rho/2)\Vert w-x+y_2/\rho\Vert^2_2-\Vert y_1\Vert^2_2/(2\rho)-\Vert y_2\Vert^2_2/(2\rho).
\end{align*}
Due to space limit, the algorithm to solve Problem \ref{prob:problem 4} is shown in Algorithm \ref{algo:neADMM spherical constraints} in the Appendix.
All subproblems are detailed as follows:\\
\textbf{1. Update $x$.}\\
\indent The variable $x$ is updated as follows:
\begin{align}
    x^{k\!+\!1}\!\leftarrow\! \arg\min\nolimits_{x} \ell(x)\!+\!(\rho/2)\Vert w^k\!-\!x\!+\!y^k_2/\rho\Vert^2_2
    \label{eq:update x}
\end{align}
This subproblem  is convex and can be solved using the Fast Iterative Shrinkage-Thresholding Algorithm (FISTA)  \cite{beck2009fast}.\\
\textbf{2. Update $w$.}\\
\indent The variable $w$ is updated as follows:
\begin{align}
    \nonumber w^{k+1}&\leftarrow  \arg\min\nolimits_{w}\Vert w-x^{k+1}+y_2^k/\rho\Vert^2_2\\&+\Vert \Vert w\Vert^2_2-1+y^k_1/\rho\Vert^2_2
    \label{eq:update w}
\end{align}
This subproblem is non-convex, but we have figured out a closed-form solution, as shown in the following theorem.
\begin{theorem}
\label{theo:sphere constraint w subproblem}
The solution to Equation \eqref{eq:update w} is
\begin{align}
w\!=\!(x^{k+1}\!-\!y^k_2/\rho)/(2\Vert w\Vert^2_2\!-\!1+2y^k_1/\rho)
\label{eq: solution to w}
\end{align}
where $u=\Vert w\Vert_2$ is obtained uniquely from one real root of the following cubic equation.
\begin{align}
\left\vert 2 u^3-u+\!2u y_1^k/\rho\right\vert=\Vert x^{k+1}-y_2^k/\rho\Vert_2 \label{eq: norm w}
\end{align}
\end{theorem}
\noindent Due to space limit, its proof is shown in Section \ref{sec: sphere constraint theorem} in the Appendix.
 \subsection{Optimization Problems with Logical Constraints with ``max'' Operations}
 \label{sec:max operator}
 \indent In the machine learning community, many minimization problems contain the ``max" operator. For example, in the multi-instance learning problem  \cite{amores2013multiple}, each ``bag'' can contain multiple ``instances'' and the classification task is to predict the labels of both bags and their instances. Here, a conventional logic between the $i$-th bag and its instances is that ``the bag label $q_i$ is 1 when at least one of its instances'' labels is 1; otherwise, the bag label $q_i$ is 0''. This well-known rule is referred to as the ``max rule''  \cite{wang2018multi}, namely $q_i=\max\nolimits_{j=1,\cdots,n_i} t_{i,j}$, where $t_{i,j}$ is the label of the $j$-th instance of $q_i$, $n_i$ is the number of instances, and $q_i$ and $t_{i,j}$ can be generalized to real numbers  \cite{amores2013multiple}. This problem can be formulated mathematically as follows:
  \begin{problem}[Multi-instance Learning Problem]
  \begin{align*}
      &\min\nolimits_{t,q,\beta} \ell(q)+\Omega(\beta) \\
     &s.t. \ q_i=\max\nolimits_{j=1,\cdots,n_i} t_{i,j},\ t_{i,j}=X_{i,j}\beta
 \end{align*}
 \label{prob:problem 7}
 \end{problem} 
 \vspace{-0.5cm}
 where $\ell(\bullet)$ and $\Omega(\bullet)$ are the loss function and the regularization term, respectively. $\beta$ is a feature weight vector, $X_{i,j}$ and $t_{i,j}$ are the $j$-th input instance and the predicted value in the $i$-th bag, and $n_i$ is the number of instances in the $i$-th bag. Let $q=[q_1,\cdots,q_n]$, $t=[t_1,\cdots,t_n]$ and $X=[X_1,\cdots,X_n]$, where $n$ is the number of bags, $t_i=[t_{i,1},\cdots,t_{i,n_i}]$, and $X_i=[X_{i,1},\cdots,X_{i,n_i}]$.\\
 \indent The augmented Lagrangian is formulated as follows according to Equation \eqref{eq: augmented Lagrangian}:
 \begin{align*}
     &L_\rho(q,\beta,t,y_1,y_2)\\&=\ell(q)+\Omega(\beta)+y_1^T(q-\max t)+(\rho/2)\Vert q-\max t\Vert^2_2\\&+y^T_2(t-X\beta)+(\rho/2)\Vert t-X\beta\Vert^2_2\\&=\ell(q)+\Omega(\beta)+(\rho/2)\Vert q-\max t+y_1/\rho\Vert^2_2\\&+(\rho/2)\Vert t-X\beta+y_2/\rho\Vert^2_2-\Vert y_1\Vert^2_2/(2\rho)-\Vert y_2\Vert^2_2/(2\rho)
 \end{align*}
 Due to space limit, the algorithm to solve Problem \ref{prob:problem 7} is shown in Algorithm \ref{algo:neADMM max constraint} in the Appendix.
 All subproblems are shown as follows:\\
\textbf{1. Update $q$.}\\
 \indent The variable $q$ is updated as follows:
 \begin{align}
       q^{k+1}&\!\leftarrow\! \arg\min\nolimits_q \ell(q)\!+\!(\rho/2)\Vert q\!-\!\max t^k\!+\!y_1^k/\rho\Vert^2_2
     \label{eq:update q}
 \end{align}
 \indent This subproblem is convex and solved using  FISTA  \cite{beck2009fast}.\\
 \textbf{2. Update $\beta$.}\\
 \indent The variable $\beta$ is updated as follows:
 \begin{align}
      \beta^{k+1}&\leftarrow\arg\min\nolimits_{\beta} \Omega(\beta)\!+\!(\rho/2)\Vert t^k\!-\!X\beta\!+\!y_2^k/\rho\Vert^2_2
     \label{eq:update beta}
 \end{align}
  \indent This subproblem is convex and solved using FISTA  \cite{beck2009fast}.\\
 \textbf{3. Update $t$.}\\
 \indent The variable $t$ is updated as follows:
 \begin{align}
     \nonumber t^{k+1}&\leftarrow \arg\min\nolimits_{t} \Vert q^{k+1}-\max t +y^{k}_1/\rho\Vert^2_2\\&+\Vert t-X\beta^{k+1}+y_2^{k}/\rho\Vert^2_2
     \label{eq:update t}
 \end{align}
This subproblem is nonconvex and difficult to solve. Here we apply a linear search to solve it. Due to the separability of t, we have
\begin{align}
\nonumber &t_i = \arg \min\nolimits_{t_i} h(t_i)\\
&h(t_i) \!=\!\Vert q^{k\!+\!1}_i \!-\! t_{i,j^*} \!+\!y^k_{1,i}/\rho\Vert_2^2 \!+\!\Vert t_i \!-\! \varphi_i 
\Vert_2^2 \label{eq:h(t_i)}\\
\nonumber &s.t. \ t_{i,j^*} = \max\nolimits_{j=1,\cdots,n_i} t_{i,j}
\end{align}
where ${\varphi_{i,j}} = X_{i,j}^T\beta^{k+1} - {y^{k}_{2,i,j}}/\rho$ is constant.\\
\indent It is easy to find that ${t_{i,j}} = \min ({t_{i,j^*}},{\varphi_{i,j}})\leq t_{i,j^*}$. We need to consider two cases: (1) $t_{i,j}<t_{i,j^*}$, (2) $t_{i,j}=t_{i,j^*}$. This problem is therefore split into two subroutines: (1) find solutions to two cases, and (2) decide which case every instance belongs to. 
\paragraph{Subroutine 1} For case (1), we have a closed-form solution $t_{i,j} = \varphi_{i,j}$. For case (2), we define a set $C = \{ j: t_{i,j}=t_{i,j^*}\}$ and its complement $\overline{C}=n_i-C$, and  we plug it in Equation \eqref{eq:h(t_i)} to obtain: 
\begin{align*}
t_{i,j^*} = (\sum\nolimits_{j \in C} \varphi_{i,j} + q^{k+1}_i + y^{k}_{1,i}/\rho)/(\vert C \vert+ 1).
\end{align*}
\paragraph{Subroutine 2}
\indent  Solving case (2) is equivalent to minimizing $h(t_i)$ by selecting appropriate indexes for set $C$. The definition of $C$ implies that $C$ consists of indexes that have 
the largest $\varphi_{i,j}$. Otherwise, $t_{i,j}<t_{i,j^*}$ and hence $i \notin C$. Now we need to determine how many indexes $C$ should have. Let $\varphi^{'}_{i}$ be a decreasing order of $\varphi_{i}$, $\vert C\vert=c$ and $a_{i,c} = (\sum\nolimits_{j = 1}^c {{\varphi^{'}_{i,j}}} + 
{q^{k+1}_i} + {y^{k}_{1,i}}/\rho)/(c + 1)$. then we have the following theorem: 
\begin{theorem} $h(t_i)|_{t_{i,j^*}=a_{i,c}}$ increases monotonically with $c$, where $h(t_i)$ is defined in Equation \eqref{eq:h(t_i)}.
\label{theo:max t subproblem}
\end{theorem}
\indent Due to space limit, its proof is in Section \ref{sec: max constraint theorem} in the Appendix. The above theorem implies that the smallest $c$ minimizes Equation \eqref{eq:h(t_i)}. So the objective becomes 
\begin{align*}
c^*=\arg \min\nolimits_c c, \ s.t. \  a_{i,c}> \varphi^{'}_{i,c+1}
\end{align*}
The corresponding solution is:
\begin{align*}
{t_{i,j}} =min(\varphi_{i,j},a_{i,c^*})
\end{align*}
The time complexity is $O(n_i log n_i)$ because the main operation of the linear search is to sort $\varphi_i$ in decreasing order.
\section{Experiments}
\indent In this section, we assess the performance of our proposed neADMM on two applications\footnote{Our code is available at \url{https://github.com/xianggebenben/neADMM}}. The experiments were conducted on a 64-bit machine equipped with an InteLoss(R) core(TM) processor (i7-6820HQ CPU) and 16.0GB memory.
\label{sec:experiments}
\subsection{1-bit Compressive Sensing}
\subsubsection{Problem Settings}
\indent  In Problem \ref{prob:problem 4}, we set $\ell(x)=\Vert x\Vert_1+\lambda/2 \sum \min(Y\phi x,0)^2$ where $\phi\in \mathbb{R}^{M\times N}$ is a measurement operator, $Y\in \mathbb{R}^{M\times M}$ is a measurement matrix and $\lambda>0$ is a tuning parameter. Here, $N$ represents the number of signals, $M$ denotes the number of measurements and $K$ denotes the number of nonzero signals. $\lambda$ and  $\rho$ were set to 0.01 and 1, respectively, and the maximal number of iterations was set to 100. The comparison methods were the Interior Point (IP) method \cite{byrd1999interior}, the Active Set (AS) method \cite{nocedal2006numerical}, and the Sequential Quadratic Programming (SQP) method \cite{fletcher2013practical}. They were all provided by the Matlab optimization toolbox and shared the same initial points.
 \subsubsection{Performance}
\begin{figure}
\begin{minipage}{0.49\linewidth} 
\centerline{\includegraphics[width=\columnwidth]{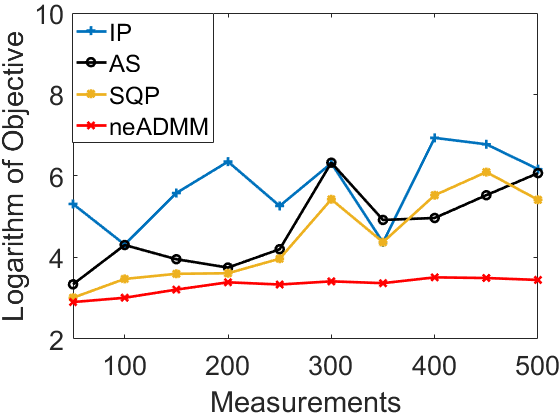}}
\centerline{(a).K=16}
\end{minipage}
\hfill
\begin{minipage}{0.49\linewidth} 
\centerline{\includegraphics[width=\columnwidth]{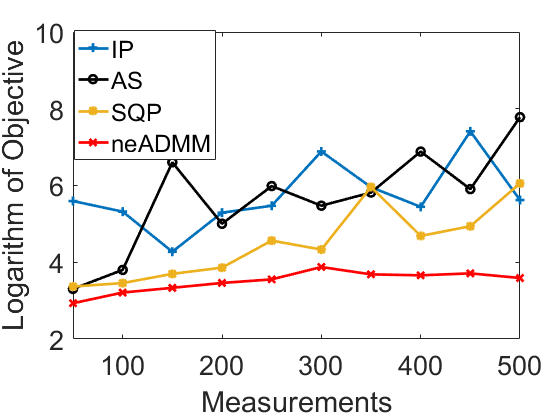}}
\centerline{(b).K=32}
\end{minipage}
\vfill
\begin{minipage}{0.49\linewidth}
\centerline{\includegraphics[width=\columnwidth]{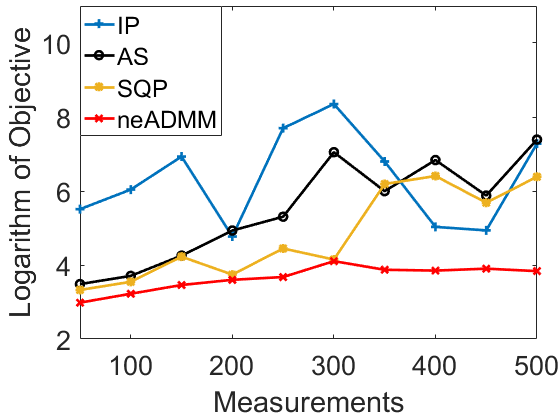}}
\centerline{(c).K=48}
\end{minipage}
\hfill
\begin{minipage}{0.49\linewidth} 
\centerline{\includegraphics[width=\columnwidth]{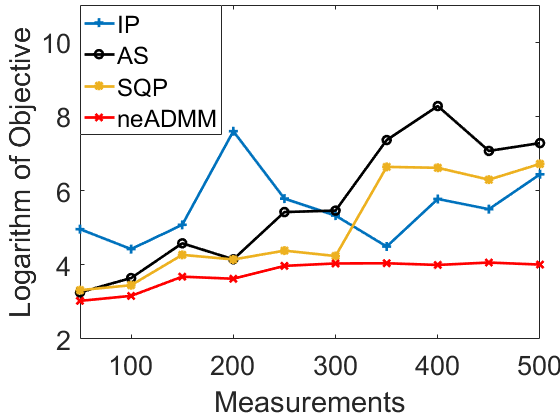}}
\centerline{(d).K=64}
\end{minipage}
\caption{Measurements versus objective values for different choices of K: the objective of the neADMM is the smallest.}
\label{fig:objective}
\vspace{-0.7cm}
\end{figure}

\indent Figure \ref{fig:objective} shows the relationship between the number of measurements and objective values for different choices of $K$. Overall, the objective values of the neADMM are lower (i.e. better) than those of the comparison methods. The logarithms of the objective values of the neADMM are all around 3, while these of comparison methods fluctuate somewhat. Even though all comparison methods are state-of-the-art, our proposed neADMM performs better maybe due to its inherent splitting schemes: all subproblems of neADMM have global optima, which make the neADMM easier to find better solutions, while all comparison methods may plunge into saddle points or local minima.\\
\subsection{Multi-instance Learning}
\indent This experiment  validates the effectiveness of our proposed neADMM against several comparison methods for multi-instance learning problems.
\subsubsection{ Problem Settings}
\indent Two datasets  are used for the performance evaluation: vaccine adverse events \cite{wang2018semi} and fox images \cite{andrews2003support}. For each dataset, $80\%$ are trained with a classifier and the remaining $20\%$ used for testing. Four comparison methods are Constructive Clustering based Ensemble (CCE) \cite{zhou2007solving}, Multi\\-instance learning with graph (miGraph) \cite{zhou2009multi}, Multi-instance Learning
based on the Vector of Locally Aggregated Descriptors representation (miVLAD) \cite{Wei2014Scalable}, and Multi-instance Learning based on the Fisher Vector representation (miFV) \cite{Wei2014Scalable}. We set the loss function $L(\bullet)$ to be a logarithm loss.  $\Omega(\beta)=\lambda \Vert\beta\Vert_1$ where $\lambda>0$ is a regularization parameter and was set to 1. The penalty parameter $\rho$ was set to 0.1, and the maximal number of iterations was set to 100. Six metrics were utilized: Accuracy (ACC), Precision (PR), Recall (RE), F-score (FS), Area Under Receiver Operating Characteristic curve (AUC) , and Area Under Precision Recall curve (AUPR). More details can be found in the Appendix.
\begin{table}
\centering
\scriptsize
\caption{Model performance on the  two datasets under six metrics.}
\begin{tabular}{p{1cm}|p{0.7cm}|p{0.7cm}|p{0.7cm}|p{0.7cm}|p{0.7cm}|p{0.7cm}}
\hline\hline
\multicolumn{7}{c}{Vaccine Adverse Events}\\\hline
Method & ACC & PR & RE & FS & AUC&AUPR \\
\hline
CCE&0.7397&0.7818&0.3805&0.5119&0.8401&0.7325\\ \hline
miGraph&0.7206&0.6812&0.4159&0.5165&0.7465&0.6128\\ \hline
miFV&0.6603&0.80&0.0708&0.1301&0.8266&0.7136 \\    \hline
miVLAD&0.7365&0.7419&0.4071&0.5257&0.7110&0.5978 \\ \hline
neADMM&\textbf{0.8000}&\textbf{0.8049}&\textbf{0.5841}&\textbf{0.6769}&\textbf{0.8901}&\textbf{0.7961}\\    \hline
\multicolumn{7}{c}{Fox Images}\\\hline
Method & ACC & PR & RE & FS & AUC&AUPR \\\hline
CCE&\textbf{0.5250}&0.4444&0.4706&0.4571&0.5358&0.4250\\ \hline
miGraph&0.4250&0.4250&\textbf{1.0000}&\textbf{0.5965}&0.4425&0.3824\\ \hline
miFV&0.4500&0.4000&0.5882&0.4762&0.4757&0.3760 \\    \hline
miVLAD&0.4500&0.4242&0.8235&0.5600
&0.5192&0.4288\\ \hline
neADMM&0.5000&\textbf{0.4483}&0.7647&0.5652&\textbf{0.5422}&\textbf{0.4583}\\\hline
\end{tabular}
\label{tab:performance}
\end{table}
\vspace{-0.3cm}
\subsubsection{Performance}
\begin{figure}
\begin{minipage}
{0.49\linewidth}
\centerline{\includegraphics[width=\columnwidth]
{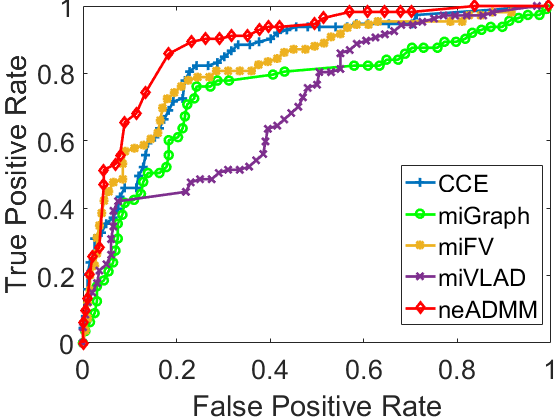}}
\centerline{(a) ROC curve}
\end{minipage}
\hfill
\begin{minipage}
{0.49\linewidth}
\centerline{\includegraphics[width=\columnwidth]
{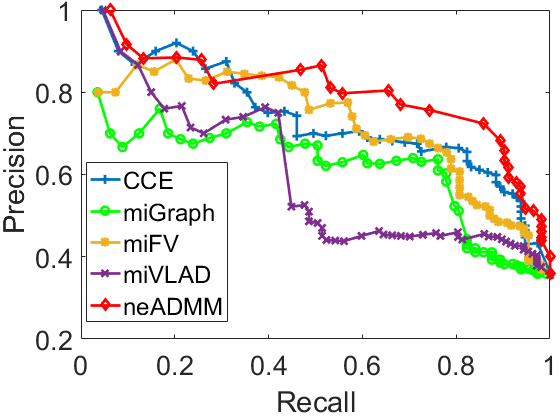}}
\centerline{(b) PR curve}
\end{minipage}
\caption{The ROC and PR curves for all methods on the vaccine adverse events dataset: the neADMM outperfoms others.}
\label{fig: vaccine roc}
\vspace{-0.5cm}
\end{figure}

\indent Table \ref{tab:performance} summarizes the prediction results obtained using neADMM and the four comparison methods on two datasets. In general, the metrics for neADMM are better than those for any of the comparison methods including for AUC and AUPR, which are the most important metrics. \\
\indent Figure \ref{fig: vaccine roc} shows that the ROC and PR curves of all the comparison methods are surrounded by these of neADMM, which is consistent with the data shown in Table \ref{tab:performance}.
\begin{figure}
\begin{minipage}
{0.49\linewidth}
\centerline{\includegraphics[width=\columnwidth]
{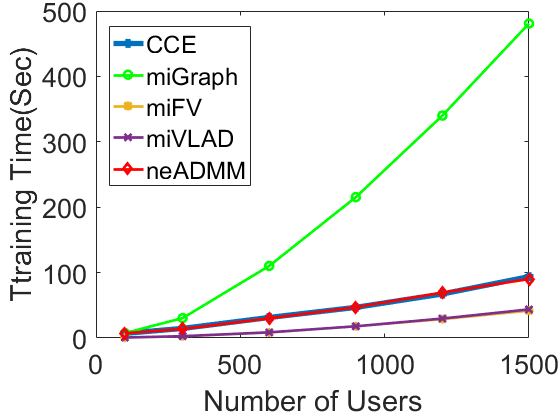}}
\centerline{(a) Time versus users}
\end{minipage}
\hfill
\begin{minipage}
{0.49\linewidth}
\centerline{\includegraphics[width=\columnwidth]
{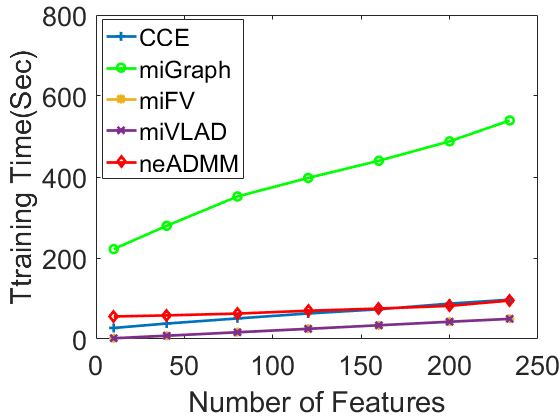}}
\centerline{(b) Time versus features}
\end{minipage}
\caption{The running time for all the methods for the vaccine adverse events dataset: the running time increases approximately linearly with the number of users and features.}
\label{fig: vaccine scalability}
\vspace{-0.7cm}
\end{figure}
\subsubsection{Scalability Analysis}
\label{sec:scalability}
\indent To examine the scalability of neADMM, we measure the running time for all methods, by averaging twenty times. As shown in Figure \ref{fig: vaccine scalability}, the running time of all the methods increases linearly with the number of users and features, and The proposed neADMM is computationally efficient.
\vspace{-0.5cm}
\section{Conclusions}
\label{sec:conclusions}
\indent We propose neADMM, an extension of ADMM framework for nonlinear equality-constrained problems. The challenge of neADMM is to solve nonconvex subproblems. Our main contribution is to provide solutions to them in two specific applications: in the spherical-constrained problem, one nonconvex subproblem is solved by roots of a cubic equation (Theorem \ref{theo:sphere constraint w subproblem}); for the problem with ``max" operations, we solve its nonconvex subproblem by a linear search (Theorem \ref{theo:max t subproblem}). We theoretically guarantee that solutions to these two subproblems are globally optimal. Experiments demonstrate that our proposed neADMM performs the best, and scales well with the increase of features and samples.\\
\indent In the future, we may explore the convergence guarantee of the proposed neADMM because convergence conditions in the existing literatures can not be applied to the neADMM.
\vspace{-0.5cm}
\bibliographystyle{abbrv}
\bibliography{Reference}
\normalsize
\newpage
\appendix
\textbf{Appendix}
\section{Theorem Proofs}
\subsection{The proof of Theorem \ref{theo:sphere constraint w subproblem}}
\label{sec: sphere constraint theorem}
\begin{proof}
\indent It is obvious that as $w\rightarrow \infty$ or $w\rightarrow -\infty$, the value of Equation \eqref{eq:update w} approaches infinity. This indicates that there exists a global minimum in Equation \eqref{eq:update w}. In general, the globally minimal point of a function is either a point 1) whose gradient is 0 or 2) whose gradient does not exist, or 3) a point at the boundary of the domain \cite{fitzpatrick2009advanced}. In our case, the domain of Equation \eqref{eq:update w} is the set of real vectors, whose boundary is empty (i.e. condition 3 is impossible), and the derivative of Equation \eqref{eq:update w} exists for any real vector $w$ (i.e. condition 2 is impossible). Therefore, the gradient of the  globally minimal point of the Equation \eqref{eq:update w} must be 0, which leads to
\begin{align}
    2(w\!-\!x^{k+1}\!+\!y^k_2/\rho)\!+\!4(\Vert w\Vert^2_2\!-\!1\!+\!2y^k_1/\rho)w=0 \label{eq: w optimality}
\end{align}
\normalsize
Equation \eqref{eq: solution to w} is obtained directly from  Equation \eqref{eq: w optimality}. Now the only issue is to obtain $\Vert w\Vert_2$, which appears on the right side of Equation \eqref{eq: solution to w}. To achieve this, we  obtain Equation \eqref{eq: norm w} by taking the  norm on both sides of Equation \eqref{eq: solution to w} as follows:
\begin{align}
    \Vert w\Vert_2=\Vert x^{k+1}\!-\!y^k_2/\rho \Vert_2/\left\rvert 2\Vert w\Vert^2_2\!-\!1+2y^k_1/\rho \right\rvert \label{eq: take norm w}
\end{align} 
where $2\Vert w\Vert^2_2\!-\!1+2y^k_1/\rho$ is a scalar.
Let $u=\Vert w\Vert_2$, Equation \eqref{eq: take norm w} is equivalent to Equation \eqref{eq: norm w}. Equation \eqref{eq: norm w} is a cubic equation with regard to $u$, and can be solved using Cardano's Formula. In order to obtain the unique $u$, we consider two possibilities of three roots of Equation \eqref{eq: norm w}: \\
1. One pair of conjugate imaginary roots and one real root. In this case, $u$ is the real root.\\
2. Three real roots. In this case, we have three possible values of $w$ using Equation \eqref{eq: solution to w}, which correspond to three real roots. $u$ is the root whose corresponding $w$  minimizes Equation \eqref{eq:update w}.\\
After $u$ is obtained, $w$ can be obtained using Equation \eqref{eq: solution to w}.
 \end{proof}
\subsection{The proof of Theorem \ref{theo:max t subproblem}}
\label{sec: max constraint theorem}
\begin{proof}
\begin{align*}
&h(t_i)|_{t_{i,j^*}=a_{i,c+1}} - h(t_i)|_{t_{i,j^*}=a_{i,c}}\\ 
&= \sum\nolimits_{j = 1}^{c + 1} {{(a_{i,c + 1} \!-\! {\varphi^{'}_{i,j}})}^2} \!+\! {(a_{i,c + 1} \!-\! {q^{k+1}_i} \!-\! {y^{k}_{1,i}}/\rho)^2}\\& - \sum\nolimits_{j = 1}^c {{{(a_{i,c} - {\varphi^{'}_{i,j}})}^2}} - {(a_{i,c} - {q^{k+1}_i} - {y^{k}_{1,i}}/\rho)^2} \\
&=(a_{i,c+1}-a_{i,c})((c+1)(a_{i,c+1}+a_{i,c})\\&-2\sum\nolimits_{j=1}^{c}\varphi^{'}_{i,j}-2q^{k+1}_i-2y^{k}_{1,i}/\rho)\\
&=(a_{i,c+1}-a_{i,c})(\varphi^{'}_{i,c+1}-a_{i,c+1})\\
&=(a_{i,c+1}-a_{i,t})((c+2)a_{i,c+1}-(c+1)a_{i,c}-a_{i,c+1})\\
  &= (c + 1){(a_{i,c + 1} - a_{i,c})^2} \geqslant 0. 
\end{align*}
Hence the theorem is proven.
 \end{proof}
 \section{Algorithms to solve Problems \ref{prob:problem 4} and \ref{prob:problem 7}}
\indent The neADMM algorithm used to solve Problem \ref{prob:problem 4} is outlined in Algorithm \ref{algo:neADMM spherical constraints}. Specifically,  Lines 9-10 update the dual variables $y_1$ and $y_2$, respectively, Lines 11-12 compute the primal and dual residuals, respectively, Lines 3-4 update $x$ and $w$ alternately. Unfortunately, Algorithm \ref{algo:neADMM spherical constraints} is not necessarily convergent by our numeric experiments, but it outperforms existing state-of-the-art methods. The computational cost of Algorithm \ref{algo:neADMM spherical constraints} mainly consists in Equation \eqref{eq:update x}, whose time complexity is $O(1/k^2)$, where $k$ is the number of iterations in FISTA \cite{beck2009fast}.  \\
 \begin{algorithm} 
  \scriptsize
\caption{The Solution to Problem \ref{prob:problem 4} Using neADMM} 
\begin{algorithmic}[1] 
\STATE Initialize $x$, $w$, $y_1$, $y_2$  $\rho>0$, $k=0$.
\REPEAT
\STATE Update $x^{k+1}$ in Equation \eqref{eq:update x}.
\STATE Update $w^{k+1}$ in Equation \eqref{eq:update w}.
\STATE Update $r_1^{k+1}\leftarrow \Vert w^{k+1}\Vert^2_2-1$.
\STATE Update $r_2^{k+1}\leftarrow  w^{k+1}-x^{k+1}$.
\STATE Update $s_1^{k+1}\leftarrow \rho(\Vert w^{k+1}\Vert^2_2-\Vert w^{k}\Vert^2_2)$.
\STATE Update $s_2^{k+1}\leftarrow \rho(w^{k+1}-w^{k})$.
\STATE Update $y_1^{k+1}\leftarrow y_1^k+\rho r_1^{k+1}$.
\STATE Update $y_2^{k+1}\leftarrow y_2^k+\rho r_2^{k+1}$.
\STATE Update $r^{k+1}\leftarrow \sqrt{\Vert r_1^{k+1}\Vert^2_2+\Vert r_2^{k+1}\Vert^2_2}$. $\#$ Calculate the primal residual.
\STATE Update $s^{k+1}\leftarrow \sqrt{\Vert s_1^{k+1}\Vert^2_2+\Vert s_2^{k+1}\Vert^2_2}$. $\#$ Calculate the dual residual.
\STATE $\nonumber k\leftarrow k+1$.
\UNTIL convergence.
\STATE Output $x$ and $w$.
\end{algorithmic}
\label{algo:neADMM spherical constraints}
\normalsize
\end{algorithm}
  \begin{algorithm}  
  \scriptsize
\caption{The Solution to Problem \ref{prob:problem 7} Using neADMM} 
\begin{algorithmic}[1] 
\STATE Initialize $q$, $\beta$, $t$, $y_1$, $y_2$,  $\rho>0$, $k=0$.
\REPEAT
\STATE Update $q^{k+1}$ in Equation \eqref{eq:update q}.
\STATE Update $\beta^{k+1}$ in Equation \eqref{eq:update beta}.
\STATE Update $t^{k+1}$ in Equation \eqref{eq:update t}.
\STATE Update $r_1^{k+1}\leftarrow q^{k+1}-\max t^{k+1}$. 
\STATE Update $r_2^{k+1}\leftarrow  t^{k+1}-X\beta^{k+1}$. 
\STATE Update $s_1^{k+1}\leftarrow \rho(\max t^k-\max t^{k+1})$. 
\STATE Update $s_2^{k+1}\leftarrow  t^{k+1}-t^k$. 
\STATE Update $y_1^{k+1}\leftarrow y_1^k+\rho r_1^{k+1}$.
\STATE Update $y_2^{k+1}\leftarrow y_2^k+\rho r_2^{k+1}$.
\STATE Update $r^{k+1}\leftarrow \sqrt{\Vert r_1^{k+1}\Vert^2_2+\Vert r_2^{k+1}\Vert^2_2}$. $\#$ Calculate the primal residual.
\STATE Update $s^{k+1}\leftarrow \sqrt{\Vert s_1^{k+1}\Vert^2_2+\Vert s_2^{k+1}\Vert^2_2}$. $\#$ Calculate the dual residual.
\STATE $\nonumber k\leftarrow k+1$.
\UNTIL convergence.
\STATE Output  $q$, $\beta$ and $t$.
\end{algorithmic}
\label{algo:neADMM max constraint}
\normalsize
\end{algorithm}
\indent The neADMM Algorithm to solve Problem \ref{prob:problem 7} is stated in Algorithm \ref{algo:neADMM max constraint}. Specifically, Lines 12-13 compute the residuals, Lines 10 and 11 update the dual variables $y_1$ and $y_2$, respectively, and Lines 3-5 update $q$,$\beta$ and $t$ alternately.
\section{More Experimental Details of  Multi-instance Learning}
\subsection{Comparison Methods}
\label{sec:comparison method}
\indent The following methods serve as baselines for the performance comparison.\\
\indent 1. Constructive Clustering-based Ensemble (CCE) \cite{zhou2007solving}. Each instance in the bag is first clustered into groups, then a classifier distinguishes a bag from the others based on group information.  Many classifiers are generated due to the many different group numbers. The final step is then to gather all the classifiers together. \\
\indent 2.   Multi-instance learning with graph (miGraph) \cite{zhou2009multi}. The miGraph treats the instances in the bags as non-independently and identically distributed.  It then implicitly constructs a graph by considering the affinity matrices and defines a new graph kernel which contains the clique information.\\
\indent 3. Multi-instance Learning
based on the Vector of Locally Aggregated Descriptors representation (miVLAD) \cite{Wei2014Scalable}. Here, multiple instances are mapped into a high dimensional vector by the Vector of Locally Aggregated Descriptors (VLAD) representation. The SVM can then be applied to train a classifier.\\
\indent 4. Multi-instance Learning
based on the Fisher Vector representation (miFV) \cite{Wei2014Scalable}. The miFV was similar to the miVLAD except that multiple instances were encoded by the Fisher Vector (FV) representation.
\subsection{Metrics}
\label{sec:metrics}
\indent This experiment utilizes six metrics to evaluate model performance: the Accuracy (ACC) is the ratio of accurately labeled bags to all bags; the Precision (PR) is the ratio of those accurately labeled as positive bags to all those labeled as positive bags; the Recall (RE) defines the ratio of those accurately labeled as positive bags to all positive bags; the F-score (FS) is the harmonic mean of the precision and recall; the Receiver Operating Characteristic curve (ROC curve) and the Precision-Recall curve (PR curve) both delineate the classification ability of a model as its discrimination threshold varies; the Area Under ROC (AUC) and the Area Under PR curve (AUPR) are the most important metrics when evaluating the performance of a classifier.

\end{document}